\documentclass{article}
\usepackage[latin5]{inputenc}
\usepackage{amssymb}
\usepackage{amsmath}
\title{ ON SOME CHARACTERIZATIONS OF RULED SURFACE OF A CLOSED TIMELIKE CURVE IN DUAL LORENTZIAN SPACE}
\numberwithin{equation}{section}

\author {Özcan BEKTAŞ $^{*}$
\and Süleyman ŞENYURT \thanks {Ordu University, Faculty of Art and Science, Department of Mathematics, 52750, Perşembe, Ordu, Turkey, ozcanbektas1986@hotmail.com, senyurtsuleyman@hotmail.com.}  }

\date{}
\begin{document}

\maketitle

\begin{abstract} In this paper, we investigate the relations between the pitch, the angle of pitch and drall of parallel ruled surface of a closed curve in dual Lorentzian space.
\end{abstract}
\textbf{Keywords:} Timelike dual curve; ruled surface; Lorentzian space; dual numbers.\\
{\it 2000 AMS Subject Classification:} 45F10, 53A04

\section{\textbf{Introduction}}

\indent  Dual numbers were introduced by W.K. Clifford (1849-79) as a tool for his geometrical investigations. After him E.Study used dual numbers and
dual vectors in his research on the geometry of lines and kinematics
[4]. The pitches and the angles of the pitches of the closed ruled
surfaces corresponds to the one parameter dual unit spherical curves in
space of lines $IR^{3} $ were calculated respectively by Hacısalihoglu
[7] and Gursoy [5].Definition of the parallel ruled surface were
presented by Blaschke (translated by Erim [3]). The integral invariants
of the parallel ruled surfaces in the 3-dimensional Euclidean space
$E^{3} $ corresponding to the unit dual spherical parallel curves were
calculated by Senyurt [11].  The set $ID=\{{\hat{\lambda}}=\lambda +\varepsilon \lambda^{*}\left| \lambda, \lambda^{*} \in IR, {\varepsilon}^{2}=0 \right. \}$ is called dual numbers set [2].\\
\indent On this set product and addition operations are respectively
\[\left(\lambda +\varepsilon \lambda ^{*} \right)+\left(\beta +\varepsilon \beta ^{*} \right)=\left(\lambda +\beta \right)+\varepsilon \left(\lambda ^{*} +\beta ^{*} \right)\]
\noindent and
\[\left(\lambda +\varepsilon \lambda ^{*} \right)\left(\beta +\varepsilon \beta ^{*} \right)=\lambda \beta +\varepsilon \left(\lambda \beta ^{*} +\lambda ^{*} \beta \right).\]

 $ID^{3} =\left\{\overrightarrow{A}=\overrightarrow{a}+\varepsilon \overrightarrow{a}^{*} \left|\overrightarrow{a},\overrightarrow{a}^{*} \in IR^{3} \right. \right\}$ the elements of $ID^{3} $are called dual vectors . On this set addition and scalar product operations are respectively\\
\indent ${\oplus :ID^{3} \times ID^{3} \to ID^{3} } \\ {\, \, \, \, \, \, \, \, \, \, \left(\overrightarrow{A},\overrightarrow{B}\right){\rm \; \; }\to \overrightarrow{A}\oplus \overrightarrow{B}=\overrightarrow{a}+\overrightarrow{b}+\varepsilon \left(\overrightarrow{a}^{*} +\overrightarrow{b}^{*} \right)}$

\[\begin{array}{l} {\odot :ID\times ID^{3} \to ID^{3} } \\ {\, \, \, \, \, \, \, \, \, \, \left(\lambda ,\overrightarrow{A}\right){\rm \; \; }\to \lambda \odot \overrightarrow{A}=\left(\lambda +\varepsilon \lambda ^{*} \right)\odot \left(\overrightarrow{a}+\varepsilon \overrightarrow{a}^{*} \right)=\lambda \overrightarrow{a}+\varepsilon \left(\lambda \overrightarrow{a}^{*} +\lambda ^{*} \overrightarrow{a}\right)} \end{array}\]
The set $\left(ID^{3} ,\oplus \right)$ is a module over the ring $\left(ID,+,\cdot \right)$. $\left(ID-Modul\right)$.

\noindent The Lorentzian inner product of dual vectors  $\overrightarrow{A}\, ,\, \overrightarrow{{\rm \; }B}\in ID^{3} \, $ is defined by

\[\left\langle \overrightarrow{A},\overrightarrow{B}\right\rangle =\left\langle \overrightarrow{a},\overrightarrow{b}\right\rangle +\varepsilon \left(\left\langle \overrightarrow{a},\overrightarrow{b}^{*} \right\rangle +\left\langle \overrightarrow{a}^{*} ,\overrightarrow{b}\right\rangle \right)\]
with the Lorentzian inner product $\overrightarrow{a}=\left(a_{1} ,a_{2} ,a_{3} \right)$ and   $\overrightarrow{b}=\left(b_{1} ,b_{2} ,b_{3} \right)\in IR^{3} $

\[\left\langle \overrightarrow{a},\overrightarrow{b}\right\rangle =-a_{1} b_{1} +a_{2} b_{2} +a_{3} b_{3} .\]
Therefore, $ID^{3} $ with the Lorentzian inner product $\left\langle \overrightarrow{A},\overrightarrow{B}\right\rangle $ is called 3-dimensional dual Lorentzian space and denoted by of  $ID_{1}^{3} =\left\{\overrightarrow{A}=\overrightarrow{a}+\varepsilon \overrightarrow{a}^{*} \left|\overrightarrow{a},\overrightarrow{a}^{*} \in IR_{1}^{3} \right. \right\}$ [14].

\indent A dual vector $\overrightarrow{A}=\overrightarrow{a}+\varepsilon \overrightarrow{a}^{*} $ $\in $ $ID_{1}^{3} $ is called A dual space-like vector  if $\left\langle \overrightarrow{A},\overrightarrow{A}\right\rangle >0$ or $\overrightarrow{A}=0$, A dual time-like vector if$\left\langle \overrightarrow{A},\overrightarrow{A}\right\rangle <0$ ,  A dual null(light-like) vector if$\left\langle \overrightarrow{A},\overrightarrow{A}\right\rangle =0$ for $\overrightarrow{A}\ne 0$ . For $\overrightarrow{A}\ne 0$, the norm $\left\| \overrightarrow{A}\right\| $ of $\overrightarrow{A}=\overrightarrow{a}+\varepsilon \overrightarrow{a}^{*} $ $\in $ $ID^{3} $ is defined by

\[\left\| \overrightarrow{A}\right\| =\sqrt{\left|\left\langle \overrightarrow{A},\overrightarrow{A}\right\rangle \right|} =\left\| \overrightarrow{a}\right\| +\varepsilon \frac{\left\langle \overrightarrow{a},\overrightarrow{a}^{*} \right\rangle }{\left\| \overrightarrow{a}\right\| } {\rm \; },\, {\rm \; }\left\| \overrightarrow{a}\right\| \ne 0 .\]
The dual Lorentzian cross-product of $\overrightarrow{A}\, ,\, \overrightarrow{{\rm \; }B}\in ID^{3} \, $is defined as

\[\overrightarrow{A}\wedge \overrightarrow{B}=\overrightarrow{a}\wedge \overrightarrow{b}+\varepsilon \left(\overrightarrow{a}\wedge \overrightarrow{b}^{*} +\overrightarrow{a}^{*} \wedge \overrightarrow{b}\right)\]
 with the Lorentzian cross-product $\overrightarrow{a}{\rm \; ,\; }\overrightarrow{b}\in IR^{3} $

\[\overrightarrow{a}\wedge \overrightarrow{b}=\left(a_{3} b_{2} -a_{2} b_{3} ,a_{1} b_{3} -a_{3} b_{1} ,a_{1} b_{2} -a_{2} b_{1} \right) [14].\]
\noindent \textbf{Theorem(E. Study):} The oriented lines in  are in one to one correspondence with the points of the dual unit sphere   where ${\rm ID}$-Modul , see [9].\\
\indent Dual number \textbf{$\Phi =\varphi +\varepsilon \varphi ^{*} $ } is called dual angle between \textbf{ $\overrightarrow{A}\, \, {\rm ve}\, \, \overrightarrow{B}\, \, $}unit dual vectors. In this place\\
\[\sinh \left(\varphi +\varepsilon \varphi ^{*} \right)=\sinh \varphi +\varepsilon \varphi ^{*} \cosh \varphi{\rm \;}\\
{\rm \; and}{\rm \;}\cosh \left(\varphi +\varepsilon \varphi^{*} \right)=\cosh \varphi +\varepsilon \varphi ^{*} \sinh \varphi. \]

\section{\textbf{ON SOME CHARACTERIZATIONS OF RULED SURFACE OF A CLOSED TIMELIKE CURVE IN DUAL LORENTZIAN SPACE $\left(ID_{1}^{3} \right)$}}

$\overrightarrow{U}=\overrightarrow{U_{1} }(t){\rm \; \; ,\; \; }\left\| \overrightarrow{U}\left(t\right)\right\| =1$ is a differentiable timelike curve in the one parameter dual unit spherical motion $K/K'$. The closed ruled surface $(\overrightarrow{U})$ corresponds to this curve in $IR^{3} $.

\noindent Let the dual orthonormal system of curve $\overrightarrow{U}=\overrightarrow{U_{1} }(t)$ as

\[\overrightarrow{U}_{1} =\overrightarrow{U}(t){\rm \; \; \; ,\; \; \; }\overrightarrow{U}_{2} =\frac{\overrightarrow{U}^{{'} } (t)}{\left\| \overrightarrow{U}(t)\right\| } {\rm \; \; \; ,\; \; \; }\overrightarrow{U}_{3} =\overrightarrow{U}_{1} \wedge \overrightarrow{U}_{2} \]

Let \textbf{$\overrightarrow{U}(t)$}be a\textbf{ }closed timelike curve with curvature $\kappa =k_{1} +\varepsilon k^{*} _{1} $ and torsion $\tau =k_{2} +\varepsilon k_{2}^{*} $ . Let Frenet frames of \textbf{$\overrightarrow{U}(t)$ }be $\left\{\overrightarrow{U_{1} },\overrightarrow{U_{2} },\overrightarrow{U_{3} }\right\}$. In this trihedron, $\overrightarrow{U_{1} }$ is timelike vector , $\overrightarrow{U_{2} }\, \, {\rm and}\, \, \overrightarrow{U_{3} }$ are spacelike vectors. For this vectors, we can write
\begin{equation} \label{GrindEQ__2_1_}
\overrightarrow{U_{1} }\wedge \overrightarrow{U_{2} }=-\overrightarrow{U_{3} }\, {\rm \; }\, \, ,\, {\rm \; }\, \, \overrightarrow{U_{2} }\wedge \overrightarrow{U_{3} }\, =\overrightarrow{U_{1} }{\rm \; \; }\, \, ,\, {\rm \; \; }\, \overrightarrow{U_{3} }\, \wedge \overrightarrow{U_{1} }\, =-\overrightarrow{U_{2} }
\end{equation}
where $\wedge $ is the Lorentzian cross product, in space \textbf{$ID_{1}^{3} $. }In this situation, the Frenet formulas are given by

\begin{equation} \label{GrindEQ__2_2_}
\overrightarrow{U_{1} }^{{'} } =\kappa \overrightarrow{U_{2} }\, \, \, \, \, ,\, \, \, \, \, \, \overrightarrow{U_{2} }^{{'} } =\kappa \overrightarrow{U_{1} }-\tau \overrightarrow{U_{3} }\, \, \, \, \, \, ,\, \, \, \, \, \overrightarrow{U_{3} }^{{'} } =\tau \overrightarrow{U_{2} }\, \, , [15].
\end{equation}
If the last equation is separated into the dual and real part, we can obtain

\begin{equation} \label{GrindEQ__2_3_}
\left\{\begin{array}{l} {\overrightarrow{u}_{1}^{'}=k_{1}\overrightarrow{u}_{2}}\\
{\overrightarrow{u}_{2}^{'}=k_{1}\overrightarrow{u}_{1}-k_{2}\overrightarrow{u}_{3}}\\
{\overrightarrow{u}_{3}^{'}=k_{2}\overrightarrow{u}_{2}}\\
{\overrightarrow{u}_{1}^{*'}=k_{1}^{*}\overrightarrow{u}_{2}+k_{1}\overrightarrow{u}_{2}^{*}}\\
{\overrightarrow{u}_{2}^{*'}=k_{1}^{*}\overrightarrow{u}_{1}-k_{2}^{*}\overrightarrow{u}_{3}+k_{1}\overrightarrow{u}_{1}^{*}-k_{2}\overrightarrow{u}_{3}^{*}}\\
{\overrightarrow{u}_{3}^{*'}=k_{2}^{*}\overrightarrow{u}_{2}+k_{2}\overrightarrow{u}_{2}^{*}}
\end{array}\right.
\end{equation}
The Frenet instantaneous rotation vector for the timelike curve is given by
\begin{equation} \label{GrindEQ__2_4_}
\overrightarrow{\Psi }=\tau \overrightarrow{U_{1} }-\kappa \overrightarrow{U_{3} }, [13].
\end{equation}
\indent In this situation for the Steiner vector of the motion, we may write
\begin{equation} \label{GrindEQ__2_5_}
\overrightarrow{D}=\oint \overrightarrow{\Psi }
\end{equation}
\noindent or
\begin{equation} \label{GrindEQ__2_6_}
\overrightarrow{D}=\overrightarrow{U}_{1}\oint\tau dt-\overrightarrow{U}_{3}\oint\kappa dt
\end{equation}
The equation \eqref{GrindEQ__2_6_} can be written type of the dual and real part as follow
\begin{equation} \label{GrindEQ__2_7_}
\left\{\begin{array}{l}
{\overrightarrow{d}=\overrightarrow{u}_{1}\oint k_{2}dt-\overrightarrow{u}_{3}\oint k_{1}dt},\\ {\overrightarrow{d}^{*}=\overrightarrow{u}_{1}^{*}\oint k_{2}dt+\overrightarrow{u}_{1}\oint k_{2}^{*}dt-\overrightarrow{u}_{3}^{*}\oint k_{1}dt-\overrightarrow{u}_{3}\oint k_{1}^{*}dt}
\end{array}\right.
\end{equation}
Now, let is calculate the integral invariants of the closed ruled surfaces  respectively. The pitch of the closed surface  is obtained as

\[L_{u_{_{1} } } =\left\langle \overrightarrow{d},\overrightarrow{u_{1} }^{*} \right\rangle +\left\langle \overrightarrow{d}^{*} ,\overrightarrow{u_{1} }\right\rangle , \]

\begin{equation} \label{GrindEQ__2_8_}
L_{u_{_{1} } } =-\oint k_{2}^{*} dt .
\end{equation}
For the dual angle of the pitch of the closed surface , we may write

\[\Lambda _{U_{1} } =-\left\langle \overrightarrow{D},\overrightarrow{U_{1} }\right\rangle , \]
Because of the equation \eqref{GrindEQ__2_6_}  we can obtain

\begin{equation} \label{GrindEQ__2_9_}
\Lambda _{U_{1} } =\oint \tau dt .
\end{equation}
If the equation \eqref{GrindEQ__2_9_} is separated into the dual and real part, we can obtain

\begin{equation} \label{GrindEQ__2_10_}
\lambda _{u_{1} } =\oint k_{2} dt \, \, \, \, \, \, ,\, \, \, \, \, \, L_{u_{_{1} } } =-\oint k_{2}^{*} dt
\end{equation}
For the drall of the closed surface , we may write

\[P_{U_{1} } =\frac{\left\langle d\overrightarrow{u_{1} },d\overrightarrow{u_{1} }^{*} \right\rangle }{\left\langle d\overrightarrow{u_{1} },d\overrightarrow{u_{1} }\right\rangle } \]
Setting by the values of the statements $d\overrightarrow{u_{1} }{\rm \; and\; }d\overrightarrow{u_{1} }^{*} $ as the equations \eqref{GrindEQ__2_3_}  into the last equations, we get

\begin{equation} \label{GrindEQ__2_11_}
P_{U_{1} } =\frac{k_{1}^{*} }{k_{1} }
\end{equation}
The pitch of the closed surface   is obtained as

\[L_{u_{_{2} } } =\left\langle \overrightarrow{d},\overrightarrow{u_{2} }^{*} \right\rangle +\left\langle \overrightarrow{d}^{*} ,\overrightarrow{u_{2} }\right\rangle ,\]

\begin{equation} \label{GrindEQ__2_12_}
L_{u_{_{2} } } =0.
\end{equation}
For the dual angle of the pitch of the closed surface , we may write

\[\Lambda _{U_{2} } =-\left\langle \overrightarrow{D},\overrightarrow{U_{2} }\right\rangle ,\]

\begin{equation} \label{GrindEQ__2_13_}
\Lambda _{U_{2} } =0.
\end{equation}
For the drall of the closed surface , we may write

\[P_{U_{2} } =\frac{\left\langle d\overrightarrow{u_{2} },d\overrightarrow{u_{2} }^{*} \right\rangle }{\left\langle d\overrightarrow{u_{2} },d\overrightarrow{u_{2} }\right\rangle } \]
Setting by the values of the statements $d\overrightarrow{u_{2} }{\rm \; and\; }d\overrightarrow{u_{2} }^{*} $  as the equations \eqref{GrindEQ__2_3_}  into the last equations, we get

\begin{equation} \label{GrindEQ__2_14_}
P_{U_{2} } =\frac{k_{2} k_{2}^{*} -k_{1} k_{1}^{*} }{k_{2}^{2} -k_{1}^{2} }
\end{equation}
The pitch of the closed surface  is obtained as

\[L_{u_{3} } =\left\langle \overrightarrow{d},\overrightarrow{u_{3} }^{*} \right\rangle +\left\langle \overrightarrow{d}^{*} ,\overrightarrow{u_{3} }\right\rangle ,\]

\begin{equation} \label{GrindEQ__2_15_}
L_{u_{3} } =-\oint k_{1}^{*} dt
\end{equation}
For the dual angle of the pitch of the closed surface , we may write

\[\Lambda _{U_{3} } =-\left\langle \overrightarrow{D},\overrightarrow{U_{3} }\right\rangle \]
Because of the equation \eqref{GrindEQ__2_6_}  we can obtain

\begin{equation} \label{GrindEQ__2_16_}
\Lambda _{U_{3} } =\oint \kappa dt
\end{equation}
If the equation \eqref{GrindEQ__2_16_} is separated into the dual and real part, we can obtain

\begin{equation} \label{GrindEQ__2_17_}
\lambda _{u_{_{3} } } =\oint k_{1} dt \, \, \, \, \, \, ,\, \, \, \, \, \, L_{u_{3} } =-\oint k_{1}^{*} dt
\end{equation}
For the drall of the closed surface , we may write

\[P_{U_{3} } =\frac{\left\langle d\overrightarrow{u_{3} },d\overrightarrow{u_{3} }^{*} \right\rangle }{\left\langle d\overrightarrow{u_{3} },d\overrightarrow{u_{3} }\right\rangle } \]
Setting by the values of the statements $d\overrightarrow{u_{3} }{\rm \; and\; }d\overrightarrow{u_{3} }^{*} $ as the equations \eqref{GrindEQ__2_3_}  into the last equations, we get

\begin{equation} \label{GrindEQ__2_18_}
P_{U_{3} } =\frac{k_{2}^{*} }{k_{2} } .
\end{equation}
Let $\Omega \left(t\right)=\omega \left(t\right)+\varepsilon \omega ^{*} \left(t\right)$be Lorentzian timelike angle of between the instantaneous dual Pfaffion vector  $\overrightarrow{\Psi }$ and the vector $\overrightarrow{U_{3} }$.

\noindent \textbf{a)}If the instantaneous dual Pfaffion vector  $\overrightarrow{\Psi }$ is spacelike ($\left|\kappa \right|>\left|\tau \right|$)

\[\kappa =\left\| \overrightarrow{\Psi }\right\| \cosh \Omega \, \, \, \, \, \, ,\, \, \, \, \, \, \tau =\left\| \overrightarrow{\Psi }\right\| \sinh \Omega \]
On the way $\overrightarrow{C}=\overrightarrow{c}+\varepsilon \overrightarrow{c}^{*} $, unit vector  about the vector $\overrightarrow{\Psi }$ direction is

\begin{equation} \label{GrindEQ__2_19_}
\overrightarrow{C}=\sinh \Omega \, \overrightarrow{U_{1} }-\cosh \Omega \, \overrightarrow{U_{3} }
\end{equation}
If the equation \eqref{GrindEQ__2_19_} is separated into the dual and real part, we can obtain

\begin{equation} \label{GrindEQ__2_20_}
\left\{\begin{array}{l} {\overrightarrow{c}=\sinh \omega \overrightarrow{u_{1} }-\cosh \omega \overrightarrow{u_{3} }} \\ {\overrightarrow{c}^{*} =\sinh \omega \overrightarrow{u_{1} }^{*} -\cosh \omega \overrightarrow{u_{3} }^{*} +\omega ^{*} \cosh \omega \overrightarrow{u_{1} }-\omega ^{*} \sinh \omega \overrightarrow{u_{3} }} \end{array}\right.
\end{equation}
The pitch of the closed surface  is obtained as

\[L_{C} =\left\langle \overrightarrow{d},\overrightarrow{c}^{*} \right\rangle +\left\langle \overrightarrow{d}^{*} ,\overrightarrow{c}\right\rangle \]

\begin{multline}\label{GrindEQ__2_21_}
L_{C} =\cosh \omega \oint k_{1}^{*} dt-\sinh \omega \oint k_{2}^{*} dt-\\
-\omega ^{*} \left(\cosh \omega \oint k_{2} dt- \sinh \omega \oint k_{1} dt \right)
\end{multline}
If we use the equations \eqref{GrindEQ__2_10_} and \eqref{GrindEQ__2_17_} into the equation \eqref{GrindEQ__2_21_} we get
\begin{equation} \label{GrindEQ__2_22_}
L_{C} =\sinh \omega L_{u_{1} } -\cosh \omega L_{u_{3} } -\omega ^{*} \left(\cosh \omega \lambda _{u_{1} } -\sinh \omega \lambda _{u_{3} } \right)
\end{equation}
For the dual angle of the pitch of the closed ruled surface  , we may write

\[\Lambda _{C} =-\left\langle \overrightarrow{D},\overrightarrow{C}\right\rangle \]
Because of the equations \eqref{GrindEQ__2_6_} and \eqref{GrindEQ__2_19_} we can obtain\textbf{  }

\begin{equation} \label{GrindEQ__2_23_}
\Lambda _{C} =\sinh \Omega \, \oint \tau dt -\cosh \Omega \, \oint \kappa dt
\end{equation}
If we use  the equations \eqref{GrindEQ__2_9_} and \eqref{GrindEQ__2_16_}   into the last equations, we get

\begin{equation} \label{GrindEQ__2_24_}
\Lambda _{C} =\sinh \Omega \, \Lambda _{U_{1} } -\cosh \Omega \, \Lambda _{U_{3} }
\end{equation}
For the drall of the closed surface , we may write

\[P_{C} =\frac{\left\langle d\overrightarrow{c},d\overrightarrow{c}^{*} \right\rangle }{\left\langle d\overrightarrow{c},d\overrightarrow{c}\right\rangle } \]

\begin{equation} \label{GrindEQ__2_25_}
P_{C} =\frac{-\omega '\omega ^{*'} +\left(k_{1} \sinh \omega -k_{2} \cosh \omega \right)\left[\left(k_{1}^{*} -k_{2} \omega ^{*} \right)\sinh \omega +\left(k_{1} \omega ^{*} -k_{2}^{*} \right)\cosh \omega \right]}{\left(k_{1} \sinh \omega -k_{2} \cosh \omega \right)^{2} -\omega '^{2} }
\end{equation}
\textbf{b)}If the instantaneous dual Pfaffion vector  $\overrightarrow{\Psi }$ is timelike ($\left|\kappa \right|<\left|\tau \right|$)

\[\kappa =\left\| \overrightarrow{\Psi }\right\| \sinh \Omega \, \, \, \, \, \, ,\, \, \, \, \, \, \tau =\left\| \overrightarrow{\Psi }\right\| \cosh \Omega \]
On the way $\overrightarrow{C}=\overrightarrow{c}+\varepsilon \overrightarrow{c}^{*} $, unit vector  about the vector $\overrightarrow{\Psi }$ direction is

\begin{equation} \label{GrindEQ__2_26_}
\overrightarrow{C}=\cosh \Omega \overrightarrow{U_{1} }-\sinh \Omega \overrightarrow{U_{3} }
\end{equation}
If the equation \eqref{GrindEQ__2_26_} is separated into the dual and real part, we can obtain

\begin{equation} \label{GrindEQ__2_27_}
\left\{\begin{array}{l} {\overrightarrow{c}=\cosh \omega \overrightarrow{u_{1} }-\sinh \omega \overrightarrow{u_{3} }} \\ {\overrightarrow{c}^{*} =\cosh \omega \overrightarrow{u_{1} }^{*} -\sinh \omega \overrightarrow{u_{3} }^{*} +\omega ^{*} \sinh \omega \overrightarrow{u_{1} }-\omega ^{*} \cosh \omega \overrightarrow{u_{3} }} \end{array}\right.
\end{equation}
The pitch of the closed surface  is obtained as

\[L_{C} =\left\langle \overrightarrow{d},\overrightarrow{c}^{*} \right\rangle +\left\langle \overrightarrow{d}^{*} ,\overrightarrow{c}\right\rangle \]

\begin{multline} \label{GrindEQ__2_28_}
L_{C} =\sinh \omega \oint k_{1}^{*} dt-\cosh \omega \oint k_{2}^{*} dt-\\-\omega ^{*} \left(\sinh \omega \oint k_{2} dt- \cosh \omega \oint k_{1} dt \right)
\end{multline}
If we use the equations \eqref{GrindEQ__2_10_} and \eqref{GrindEQ__2_17_} into the equation \eqref{GrindEQ__2_21_} we get

\begin{equation} \label{GrindEQ__2_29_}
L_{C} =\cosh \omega L_{u_{1} } -\sinh \omega L_{u_{3} } -\omega ^{*} \left(\sinh \omega \lambda _{u_{1} } -\cosh \omega \lambda _{u_{3} } \right)
\end{equation}
For the dual angle of the pitch of the closed ruled surface  , we may write

\[\Lambda _{C} =-\left\langle \overrightarrow{D},\overrightarrow{C}\right\rangle \]
Because of the equations \eqref{GrindEQ__2_6_} and \eqref{GrindEQ__2_19_} we can obtain\textbf{ }

\begin{equation} \label{GrindEQ__2_30_}
\Lambda _{C} =\cosh \Omega \, \oint \tau dt -\sinh \Omega \, \oint \kappa dt
\end{equation}
If we use the equations \eqref{GrindEQ__2_9_} and \eqref{GrindEQ__2_16_}  into the last equations, we get

\begin{equation} \label{GrindEQ__2_31_}
\Lambda _{C} =\cosh \Omega \, \Lambda _{U_{1} } -\sinh \Omega \, \Lambda _{U_{3} }
\end{equation}
For the drall of the closed surface , we may write

\[P_{C} =\frac{\left\langle d\overrightarrow{c},d\overrightarrow{c}^{*} \right\rangle }{\left\langle d\overrightarrow{c},d\overrightarrow{c}\right\rangle } \]

\begin{equation} \label{GrindEQ__2_32_}
P_{C} =\frac{\omega '\omega ^{*'} +\left(k_{1} \cosh \omega -k_{2} \sinh \omega \right)\left[\left(k_{1}^{*} -k_{2} \omega ^{*} \right)\cosh \omega +\left(k_{1} \omega ^{*} -k_{2}^{*} \right)\sinh \omega \right]}{\left(k_{1} \sinh \omega -k_{2} \cosh \omega \right)^{2} -\omega '^{2} }
\end{equation}
\textbf{Definition:}The closed ruled surface $(\overrightarrow{U})$ corresponds to the dual timelike curve \textbf{$\overrightarrow{U}(t)$ }which makes the fixed dual angle\textbf{ }$\Phi =\varphi +\varepsilon \varphi^{*} $ with $\overrightarrow{U}(t)$ and defines by

\begin{equation} \label{GrindEQ__2_33_}
\overrightarrow{V}=\cosh \Phi \overrightarrow{U_{1} }+\sinh \Phi \overrightarrow{U_{3} }
\end{equation}
This surface ($\overrightarrow{V}$) corresponds to dual timelike vector $\overrightarrow{V}$ is called the parallel ruled surface of surface $(\overrightarrow{U})$ in dual lorentzian space \textbf{$ID_{1}^{3} $.}

Let be $\overrightarrow{V}_{1} =\overrightarrow{V}$. Differentiating of the vector $\overrightarrow{V_{1} }$ with respect the parameter and using the equation \eqref{GrindEQ__2_3_} we get

\begin{equation} \label{GrindEQ__2_34_}
\overrightarrow{V_{1} }^{{'} } =\left(\kappa \cosh \Phi +\tau \sinh \Phi \right)\overrightarrow{U_{2} }
\end{equation}
If the norm of the vector  denotes by $P$,  we get

\begin{equation} \label{GrindEQ__2_35_}
\overrightarrow{P} =\kappa \cosh \Phi +\tau \sinh \Phi
\end{equation}
Then if is known that

\noindent Substituting by the values of the equations \eqref{GrindEQ__2_34_} and \eqref{GrindEQ__2_35_} into   , we get

\begin{equation} \label{GrindEQ__2_36_}
\overrightarrow{V_{2} }=\overrightarrow{U_{2} }
\end{equation}
Then, fort he vector   , we get

\begin{equation} \label{GrindEQ__2_37_}
\overrightarrow{V_{3} }=-\sinh \Phi \overrightarrow{U_{1} }-\cosh \Phi \overrightarrow{U_{3} }
\end{equation}
If the equation  \eqref{GrindEQ__2_33_} , \eqref{GrindEQ__2_36_}  and  \eqref{GrindEQ__2_37_}   are written matrix form, we have
\begin{eqnarray}\label{GrindEQ__2_38_}
\left[\begin{array}{c}
\overrightarrow{V_{1}}\\
\overrightarrow{V_{2}}\\
\overrightarrow{V_{3}}
\end{array}
\right]
=\left[\begin{array}{ccc}
\cosh \Phi & 0 & \sinh \Phi\\
0 & 1 & 0\\
-\sinh \Phi & 0 & -\cosh \Phi
\end{array}
\right]
\cdot
\left[\begin{array}{c}
\overrightarrow{U_{1}}\\
\overrightarrow{U_{2}}\\
\overrightarrow{U_{3}}
\end{array}
\right]
\end{eqnarray}
\noindent or
\begin{eqnarray}\label{GrindEQ__2_39_}
\left[\begin{array}{c}
\overrightarrow{U_{1}}\\
\overrightarrow{U_{2}}\\
\overrightarrow{U_{3}}
\end{array}
\right]
=\left[\begin{array}{ccc}
\cosh \Phi & 0 & \sinh \Phi\\
0 & 1 & 0\\
-\sinh \Phi & 0 & -\cosh \Phi
\end{array}
\right]
\cdot
\left[\begin{array}{c}
\overrightarrow{V_{1}}\\
\overrightarrow{V_{2}}\\
\overrightarrow{V_{3}}
\end{array}
\right]
\end{eqnarray}
If the equation \eqref{GrindEQ__2_39_} is separated into real and dual parts, we get
\begin{equation} \label{GrindEQ__2_40_}
\left\{\begin{array}{l} {\overrightarrow{u}_{1}=\cosh \varphi\overrightarrow{v}_{1}+ \sinh \varphi \overrightarrow{v}_{3}}\\
{\overrightarrow{u}_{2}=\overrightarrow{v}_{2}}\\
{\overrightarrow{u}_{3}=- \sinh \varphi\overrightarrow{v}_{1}- \cosh \varphi\overrightarrow{v}_{3}}\\
{\overrightarrow{u}_{1}^{*}=\cosh \varphi\overrightarrow{v}_{1}^{*}+ \sinh \varphi \overrightarrow{v}_{3}^{*}+ \varphi^{*}( \sinh \varphi\overrightarrow{v}_{1}+ \cosh \varphi\overrightarrow{v}_{3})}\\
{\overrightarrow{u}_{2}^{*}=\overrightarrow{v}_{2}^{*}}\\
{\overrightarrow{u}_{3}^{*}=- \sinh \varphi\overrightarrow{v}_{1}^{*}- \cosh \varphi\overrightarrow{v}_{3}^{*}- \varphi^{*}( \cosh \varphi\overrightarrow{v}_{1}+ \sinh \varphi\overrightarrow{v}_{3})} \end{array}
\right.
\end{equation}
Let  be curvature $P=\, p+\varepsilon p^{*} $and torsion $Q=q+\varepsilon q^{*} $ of curve \textbf{$\overrightarrow{V}(t)$. }Between the vectors $\overrightarrow{V_{1} },\overrightarrow{V_{2} },\overrightarrow{V_{3} }$ and the derivate vectors $\overrightarrow{V_{1} }^{{'} } ,\overrightarrow{V_{2} }^{{'} } ,\overrightarrow{V_{3} }^{{'} } $ there is following relation

\begin{equation} \label{GrindEQ__2_41_}
\left\{\begin{array}{l} {{\overrightarrow{V}_{1}}^{'}=P{\overrightarrow{V}_{2}}, {\rm \; \; \; \; \; \; } {\overrightarrow{V}_{2}}^{'}=P\overrightarrow{V}_{1}-Q\overrightarrow{V}_{3},{\rm \; \; \; \; \; \; } {\overrightarrow{V}_{3}}^{'}=Q\overrightarrow{V}_{2}}\\
{P=\sqrt{<{\overrightarrow{V}_{1}}^{'},{\overrightarrow{V}_{1}}^{'}>}, {\rm \; \; \; \; \; \; } Q=\frac{det(\overrightarrow{V}_{1},{\overrightarrow{V}_{1}}^{'},{\overrightarrow{V}_{1}}^{''})}{<{\overrightarrow{V}_{1}}^{'},{\overrightarrow{V}_{1}}^{'}>}}.
\end{array}\right. [15].
\end{equation}
If the equation \eqref{GrindEQ__2_41_} is separated into the real and dual parts, we can write
\begin{equation} \label{GrindEQ__2_42_}
\left\{\begin{array}{l}{{\overrightarrow{v}_{1}}^{'}=p\overrightarrow{v}_{2}, {\rm \; \; } {\overrightarrow{v}_{2}}^{'}=p\overrightarrow{v}_{1}-q\overrightarrow{v}_{3}, {\rm \; \; } {\overrightarrow{v}_{3}}^{'}=q\overrightarrow{v}_{2}}\\
{\overrightarrow{v}_{1}^{'}{}^{*}=p\overrightarrow{v}_{2}^{*}+{p}^{*}\overrightarrow{v}_{2},} \\ {\overrightarrow{v}_{2}^{'}{}^{*}=p\overrightarrow{v}_{1}^{*}+{p}^{*}\overrightarrow{v}_{1}-q\overrightarrow{v}_{3}^{*}
-{q}^{*}\overrightarrow{v}_{3},}\\ {\overrightarrow{v}_{3}^{'}{}^{*}=q\overrightarrow{v}_{2}^{*}+{q}^{*}\overrightarrow{v}_{2}}
\end{array}\right.
\end{equation}
Now, we can calculate the value of $Q\, {\rm relative\; to\; }\kappa \, {\rm and\; }\tau $. Derivative  the equation \eqref{GrindEQ__2_34_} with respect to the parameter an of making the recesiory operations we may write
\begin{equation} \label{GrindEQ__2_43_}
\begin{split}
\overrightarrow{V}_{1}^{''}= & ({\kappa}^{2}\cosh \Phi+ \kappa\tau\sinh \Phi)\overrightarrow{U}_{1}+\\
& + (\kappa\cosh\Phi+\tau\sinh\Phi)^{'}\overrightarrow{U}_{2}+
(-\kappa\tau\cosh\Phi-{\tau}^{2}\sinh\Phi)\overrightarrow{U}_{3}
\end{split}
\end{equation}
Substituting by the equations \eqref{GrindEQ__2_33_} , \eqref{GrindEQ__2_34_} and \eqref{GrindEQ__2_43_}  into the equation \eqref{GrindEQ__2_41_}, we get
\begin{equation} \label{GrindEQ__2_44_}
Q=-\kappa\sinh\Phi-\tau\cosh\Phi
\end{equation}
The equations \eqref{GrindEQ__2_35_} and \eqref{GrindEQ__2_44_}  are separated into the dual and real parts, we have
\begin{equation} \label{GrindEQ__2_45_}
\left\{\begin{array}{l}{p{\rm\;\;}=k_{1}\cosh\varphi+k_{2}\sinh\varphi}\\
{p^{*}=k_{1}^{*}\cosh\varphi+k_{2}^{*}\sinh\varphi+{\varphi}^{*}(k_{1}\sinh\varphi+k_{2}\cosh\varphi)}\\
{q{\rm\;\;}=-k_{1}\sinh\varphi-k_{2}\cosh\varphi}\\
{q^{*}=-k_{1}^{*}\sinh\varphi-k_{2}^{*}\cosh\varphi-{\varphi}^{*}(k_{1}\cosh\varphi+k_{2}\sinh\varphi)}
\end{array}\right.
\end{equation}
In the dual unit spharical motion $K/K'$, the dual orthonormal system  $\left\{\overrightarrow{V_{1} },\overrightarrow{V_{2} },\overrightarrow{V_{3} }\right\}$ each time t is make a dual rotation motion around the instantaneous dual Pfaffion vector. This vector is determined the following equation
\begin{equation} \label{GrindEQ__2_46_}
\overrightarrow{\overline{\Psi }}=Q\overrightarrow{V_{1} }-P\overrightarrow{V_{3} } , [13]
\end{equation}
For the Steiner vector of the motion, we can write
\begin{equation} \label{GrindEQ__2_47_}
\overrightarrow{\overline{D}}=\oint \overrightarrow{\overline{\Psi }}
\end{equation}
\noindent or
\begin{equation} \label{GrindEQ__2_48_}
\overrightarrow{\overline{D}}=\overrightarrow{V_{1} }\oint Qdt-\overrightarrow{V_{3} }\oint Pdt
\end{equation}
Setting by the values of the vectors  $\overrightarrow{U_{1} }\, {\rm and\; \; }\overrightarrow{U_{3} }$ as the equations \eqref{GrindEQ__2_40_} into the equations \eqref{GrindEQ__2_4_}, we get
\begin{equation*}
\overrightarrow{\Psi }=-Q\overrightarrow{V_{1}}+P\overrightarrow{V_{3}}
\end{equation*}
\noindent In this case, we consider the last equation and equation \eqref{GrindEQ__2_46_}, we can write
\begin{equation} \label{GrindEQ__2_49_}
\overrightarrow{\Psi }=-\overrightarrow{\overline{\Psi }}
\end{equation}
Because of the equations $\overrightarrow{D}=\oint \overrightarrow{\Psi } $ and $\overrightarrow{\overline{D}}=\oint \overrightarrow{\overline{\Psi }} $ we can obtain $\overrightarrow{D}=-\overrightarrow{\overline{D}}$ . Then, for the dual Steiner vector of the motion, we may write
\begin{equation} \label{GrindEQ__2_50_}
\overrightarrow{D}=-\overrightarrow{V_{1} }\oint Qdt+\overrightarrow{V_{3} }\oint Pdt
\end{equation}
It is separated into the dual and real part as
\begin{equation} \label{GrindEQ__2_51_}
\left\{\begin{array}{l}\overrightarrow{d}=-\overrightarrow{v}_{1}\oint qdt+\overrightarrow{v}_{3}\oint pdt,\\ \overrightarrow{d}^{*}=-\overrightarrow{v}_{1}\oint q^{*}dt-\overrightarrow{v}_{1}^{*}\oint qdt+\overrightarrow{v}_{3}\oint {p}^{*}dt+\overrightarrow{v}_{3}^{*}\oint pdt
\end{array}\right.
\end{equation}
Now, let is calculate the integral invariants of the closed ruled surfaces  respectively. The pitch of the closed surface  is obtained as
\[L_{V_{1} } =\left\langle \overrightarrow{d},\overrightarrow{v_{1} }^{*} \right\rangle +\left\langle \overrightarrow{d}^{*} ,\overrightarrow{v_{1} }\right\rangle ,\]
\begin{equation} \label{GrindEQ__2_52_}
L_{V_{1} } =\oint q^{*} dt .
\end{equation}
Substituting by the value  into the equation \eqref{GrindEQ__2_52_}
\begin{multline} \label{GrindEQ__2_53_}
L_{V_{1}}=-\sinh \varphi\oint k_{1}^{*}dt-\cosh\varphi\oint k_{2}^{*}dt-\\
-{\varphi}^{*}\left(\cosh \varphi\oint k_{1}dt+\sinh\varphi\oint k_{2}dt \right).
\end{multline}
\noindent or
\begin{equation} \label{GrindEQ__2_54_}
L_{V_{1} } =\cosh \varphi L_{u_{_{1} } } +\sinh \varphi L_{u_{_{3} } } -\varphi ^{*} \left(\sinh \varphi \lambda _{u_{_{1} } } +\cosh \varphi \lambda _{u_{3}} \right).
\end{equation}
For the dual angle of the pitch of the closed ruled surface  , we may write
\[\Lambda _{V_{1} } =-\left\langle \overrightarrow{D},\overrightarrow{V_{1} }\right\rangle \]
Because of the equation \eqref{GrindEQ__2_50_} we can obtain
\begin{equation} \label{GrindEQ__2_55_}
\Lambda _{V_{1} } =-\oint Qdt .
\end{equation}
Substituting by the equation \eqref{GrindEQ__2_44_} into the last equation, we get
\[\Lambda _{V_{1} } \, =\sinh \Phi \oint \kappa  dt+\cosh \Phi \oint \tau  dt\]
\noindent or
\begin{equation} \label{GrindEQ__2_56_}
{\wedge}_{V_{1}}=\cosh\Phi{\wedge}_{U_{1}}+\sinh\Phi{\wedge}_{U_{3}}
\end{equation}
If the equation \eqref{GrindEQ__2_56_} is separated into the dual and real part, we can obtain
\begin{equation} \label{GrindEQ__2_57_}
\left\{\begin{array}{l} {\lambda _{v_{1} } =\cosh \varphi \lambda _{u_{_{1} } } +\sinh \varphi \lambda _{u_{_{3} } } \, \, \, \, \, \, } \\ {L_{v_{1} } =\cosh \varphi L_{u_{_{1} } } +\sinh \varphi L_{u_{_{3} } } -\varphi ^{*} \left(\sinh \varphi \lambda_{u_{1}} +\cosh \varphi \lambda_{u_{3}}\right)} \end{array}\right.
\end{equation}
For the drall of the closed surface , we may write
\[P_{V_{1} } =\frac{\left\langle d\overrightarrow{v_{1} },d\overrightarrow{v_{1} }^{*} \right\rangle }{\left\langle d\overrightarrow{v_{1} },d\overrightarrow{v_{1} }\right\rangle } \]
Setting by the values of the statements $d\overrightarrow{v_{1} }{\rm \; and\; }d\overrightarrow{v_{1} }^{*} $ as the equations \eqref{GrindEQ__2_42_}  into the last equations, we get
\begin{equation} \label{GrindEQ__2_58_}
P_{V_{1} } =\frac{p^{*} }{p}
\end{equation}
Setting by the values of $p\, {\rm and\; p}^{*} $  as the equations \eqref{GrindEQ__2_45_}  into the last equations, we get
\begin{equation} \label{GrindEQ__2_59_}
P_{V_{1} } =\frac{k_{1}^{*} \cosh \varphi +k_{2}^{*} \sinh \varphi }{k_{1} \cosh \varphi+k_{2} \sinh \varphi}+{\varphi} ^{*} \frac{k_{1} \sinh \varphi +k_{2} \cosh \varphi}{k_{1} \cosh \varphi +k_{2} \sinh \varphi}
\end{equation}
\textbf{Theorem 2.1: }Let $\left(V_{1} \right)$ be the parallel surface of the surface $\left(U_{1} \right)$. {\rm\;\;\;}The pitch, drall and the dual of the pitch of the ruled surface $\left(V_{1} \right)$ are
\[ 1-) L_{V_{1}}=\oint q^{*} dt {\rm \; \; \; \; \; \; \; \; \;} 2-)\Lambda_{V_{1}}=-\oint Qdt {\rm \; \; \; \; \; \; \; \; \;}         3-)P_{V_{1}} =\frac{p^{*} }{p}. \]
\textbf{Corollary 2.1: }Let $\left(V_{1} \right)$ be the parallel surface of the surface $\left(U_{1} \right)$. The pitch and the dual of the pitch of the ruled surface $\left(V_{1} \right)$ related to the invariants of the surface $\left(U_{1} \right)$ are written as follow\\
\indent $1-) L_{V_{1} } =\cosh \varphi L_{u_{_{1} } } +\sinh \varphi L_{u_{_{3} } } -\varphi ^{*} \left(\sinh \varphi \lambda _{u_{_{1} } } +\cosh \varphi\lambda _{u_{3}} \right)$\\
\indent $2-) {\wedge}_{V_{1}}=\cosh\Phi {\wedge}_{U_{1}}+\sinh\Phi{\wedge}_{U_{3}}$\\
The pitch of the closed surface $\left(V_{2} \right)$ is obtained as
\[L_{V_{2} } =\left\langle \overrightarrow{d},\overrightarrow{v_{2} }^{*} \right\rangle +\left\langle \overrightarrow{d}^{*} ,\overrightarrow{v_{2} }\right\rangle \]
\begin{equation} \label{GrindEQ__2_60_}
L_{V_{2} } =0
\end{equation}
For the dual angle of the pitch of the closed ruled surface $\left(V_{2} \right)$, we may write
\[\Lambda _{V_{2} } =-\left\langle \overrightarrow{D},\overrightarrow{V_{2} }\right\rangle \]
Because of the equation \eqref{GrindEQ__2_32_} we can obtain
\begin{equation} \label{GrindEQ__2_61_}
\Lambda _{V_{2} } =0
\end{equation}
For the drall of the closed surface $\left(V_{2} \right)$, we may write
\[P_{V_{2} } =\frac{\left\langle d\overrightarrow{v_{2} },d\overrightarrow{v_{2} }^{*} \right\rangle }{\left\langle d\overrightarrow{v_{2} },d\overrightarrow{v_{2} }\right\rangle } \]
Setting by the values of the statements $d\overrightarrow{v_{2} }{\rm \; and\; }d\overrightarrow{v_{2} }^{*} $ as the equations \eqref{GrindEQ__2_24_}  into the last equations, we get
\begin{equation} \label{GrindEQ__2_62_}
P_{V_{2} } =\frac{qq^{*} -pp^{*} }{q^{2} -p^{2} }
\end{equation}
Setting by the values of $p,p^{*} ,q\, {\rm and\; }q^{*} $ as the equations \eqref{GrindEQ__2_45_}  into the last equations, we get
\begin{equation} \label{GrindEQ__2_63_}
P_{V_{2} } =\frac{k_{2} k_{2}^{*} -k_{1} k_{1}^{*} }{k_{2}^{2} -k_{1}^{2} }
\end{equation}
\textbf{Theorem 2.2:} Let $\left(V_{1} \right)$ be the parallel surface of the surface $\left(U_{1} \right)$. {\rm \; \; \;}The pitch, drall and the dual of the pitch of the ruled surface $\left(V_{2} \right)$ are
\[ 1-) L_{V_{2}}=0{\rm \; \; \; \; \; \; \;\;}2-)\Lambda _{V_{2} }=0{\rm \; \; \; \; \; \; \;\;}3-) P_{V_{2} } =\frac{qq^{*} -pp^{*} }{q^{2} -p^{2} } \]
The pitch of the closed surface $\left(V_{3} \right)$ is obtained as
\[L_{V_{_{3} } } =\left\langle \overrightarrow{d},\overrightarrow{v_{3} }^{*} \right\rangle +\left\langle \overrightarrow{d}^{*} ,\overrightarrow{v_{3} }\right\rangle ,\]
\begin{equation} \label{GrindEQ__2_64_}
L_{V_{_{3} } } =\oint p^{*} dt
\end{equation}
Substituting by the value  into the equation \eqref{GrindEQ__2_64_}
\begin{equation} \label{GrindEQ__2_65_}
L_{V_{3}}=\cosh\varphi\oint k_{1}^{*}dt+\sinh\varphi\oint k_{2}^{*}dt+{\varphi}^{*}(\sinh\varphi\oint k_{1}dt+\cosh\varphi\oint k_{2}dt)
\end{equation}
\noindent or
\begin{equation} \label{GrindEQ__2_66_}
L_{V_{3}} =-\sinh \varphi L_{u_{_{1} } } -\cosh \varphi L_{u_{3}} +\varphi ^{*} \left(\cosh \varphi \lambda _{u_{1}} +\sinh \varphi \lambda_{u_{3}}\right)
\end{equation}
For the dual angle of the pitch of the closed ruled surface  , we may write
\[\Lambda _{V_{3} } =-\left\langle \overrightarrow{D},\overrightarrow{V_{3} }\right\rangle \]
Because of the equation \eqref{GrindEQ__2_50_} we can obtain
\begin{equation} \label{GrindEQ__2_67_}
\Lambda _{V_{3} } =-\oint Pdt .
\end{equation}
Substituting by the equation \eqref{GrindEQ__2_35_} into the last equation, we get
\[\Lambda _{V_{3} } \, =-\cosh \Phi \oint \kappa  dt-\sinh \Phi \oint \tau  dt\]
\noindent or
\begin{equation} \label{GrindEQ__2_68_}
{\wedge}_{V_{3}}=-\sinh\Phi{\wedge}_{U_{1}}-\cosh\Phi{\wedge}_{U_{3}}
\end{equation}
If the equation \eqref{GrindEQ__2_68_} is separated into the dual and real part, we can obtain
\begin{equation} \label{GrindEQ__2_69_}
\left\{\begin{array}{l} {\lambda_{v_{3}}}=-\sinh \varphi \lambda _{u_{1}}-\cosh \varphi \lambda_{u_{3}}\\ {{L_{v_{3}}}=-\sinh \varphi L_{u_{1}}-\cosh \varphi L_{u_{3}}+{\varphi}^{*} \left(\cosh \varphi \lambda_{u_{1}}+\sinh \varphi \lambda _{u_{3}}\right)} \end{array}\right.
\end{equation}
For the drall of the closed surface , we may write
\[P_{V_{3} } =\frac{\left\langle d\overrightarrow{v_{3} },d\overrightarrow{v_{3} }^{*} \right\rangle }{\left\langle d\overrightarrow{v_{3} },d\overrightarrow{v_{3} }\right\rangle } \]
Setting by the values of the statements $d\overrightarrow{v_{3} }{\rm \; and\; }d\overrightarrow{v_{3} }^{*} $ as the equations \eqref{GrindEQ__2_42_}  into the last equations, we get
\begin{equation} \label{GrindEQ__2_70_}
P_{V_{3} } =\frac{q^{*} }{q}
\end{equation}
Setting by the values of $q\, {\rm and\; q}^{*} $  as the equations \eqref{GrindEQ__2_45_}  into the last equations, we get
\begin{equation} \label{GrindEQ__2_71_}
P_{V_{3} } =\frac{-k_{1}^{*} \sinh \varphi -k_{2}^{*} \cosh \varphi}{-k_{1} \sinh \varphi -k_{2} \cosh \varphi }-{\varphi} ^{*} \left(\frac{k_{1} \cosh \varphi +k_{2} \sinh \varphi }{-k_{1} \sinh \varphi -k_{2} \cosh \varphi } \right)
\end{equation}
\textbf{Theorem 2.3: }Let $\left(V_{1} \right)$ be the parallel surface of the surface $\left(U_{1} \right)$. The pitch , drall and the dual of the pitch of the ruled surface $\left(V_{3} \right)$ are
\[ 1-) L_{V_{3}}=\oint p^{*}dt{\rm \; \; \;\; \; \;\; \; \;} 2-)  \Lambda _{V_{3}}=-\oint Pdt{\rm \; \; \;\; \; \;\; \; \;} 3-)  P_{V_{3}}=\frac{q^{*}}{q} \]
\textbf{Corollary 2.2: }Let $\left(V_{1} \right)$ be the parallel surface of the surface $\left(U_{1} \right)$. The pitch and the dual of the pitch of the ruled surface $(V_{3})$ related to the invariants of the surface $(U_{1})$ are written as follow\\
\indent $1-)  L_{V_{3}} =-\sinh \varphi L_{u_{1}} -\cosh \varphi L_{u_{3}}+\varphi ^{*}(\cosh \varphi \lambda _{u_{1}} +\sinh \varphi \lambda _{u_{3}} )$\\
\indent $2-) {\wedge}_{V_{3}}=-\sinh\Phi{\wedge}_{U_{1}}-\cosh\Phi{\wedge}_{U_{3}}$\\
\noindent Let $\Theta(t)=\theta(t)+\varepsilon \theta ^{*} (t)$ be Lorentzian timelike angle of between the instantaneous dual Pfaffion vector  $\overrightarrow{\overline{\Psi }}$ and the vector $\overrightarrow{V_{3}}$.

\noindent \textbf{a)}If the instantaneous dual Pfaffion vector  $\overrightarrow{\overline{\Psi }}$ is spacelike ($\left|P\right|>\left|Q\right|$)

\[P=\left\| \overrightarrow{\overline{\Psi }}\right\| \cosh \Theta \, \, \, \, \, \, ,\, \, \, \, \, \, \, \, Q=\left\| \overrightarrow{\overline{\Psi }}\right\| \sinh \Theta \]
On the way $\overrightarrow{\overline{C}}=\overrightarrow{\overline{c}}+\varepsilon \overrightarrow{\overline{c}}^{*} $, unit vector  about the vector $\overrightarrow{\overline{\Psi }}$ direction is

\begin{equation} \label{GrindEQ__2_72_}
\overrightarrow{\overline{C}}=\sinh \Theta \, \overrightarrow{V_{1} }-\cosh \Theta \, \overrightarrow{V_{3} }
\end{equation}
Setting by the values of the vectors  $\overrightarrow{V_{1} }{\rm \; and\; }\overrightarrow{V_{3} }$ as the equations \eqref{GrindEQ__2_38_}  into the equation \eqref{GrindEQ__2_72_}, we get

\[\overrightarrow{\overline{C}}=\left(\sinh \Theta \cosh \Phi +\cosh \Theta \sinh \Phi \right)\overrightarrow{U_{1} }+\left(\cosh \Theta \cosh \Phi +\sinh \Theta \sinh \Phi \right)\overrightarrow{U_{3} }\]

\begin{equation} \label{GrindEQ__2_73_}
\overrightarrow{\overline{C}}=\sinh \left(\Theta +\Phi \right)\overrightarrow{U_{1} }+\cosh \left(\Theta +\Phi \right)\overrightarrow{U_{3} }
\end{equation}
If the equation \eqref{GrindEQ__2_72_} is separated into the dual and real part, we can obtain

\begin{equation} \label{GrindEQ__2_74_}
\left\{\begin{array}{l} {\overrightarrow{\overline{c}}=\sinh \theta \overrightarrow{v_{1} }-\cosh \theta \overrightarrow{v_{3} }} \\ {\overrightarrow{\overline{c}}^{*} =\sinh \theta \overrightarrow{v_{1} }^{*} -\cosh \theta \overrightarrow{v_{3} }^{*} +\theta ^{*} \cosh \theta \overrightarrow{v_{1} }-\theta ^{*} \sinh \theta \overrightarrow{v_{3} }} \end{array}\right.
\end{equation}
The pitch of the closed surface  is obtained as

\[L_{\overline{C}} =\left\langle \overrightarrow{d},\overrightarrow{\overline{c}}^{*} \right\rangle +\left\langle \overrightarrow{d}^{*} ,\overrightarrow{\overline{c}}\right\rangle \]
Setting by the values of the statements $\overrightarrow{d}{\rm \; and\; }\overrightarrow{d}^{*} $ as the equations \eqref{GrindEQ__2_51_}  into the last equations and if we do the necessary operations , we get

\begin{equation} \label{GrindEQ__2_75_}
L_{\overline{C}} =-\cosh \theta \oint p^{*} dt +\sinh \theta \oint q^{*} dt +\theta ^{*} \left(\cosh \theta \oint qdt- \sinh \theta \oint pdt \right)
\end{equation}
If we use the equations \eqref{GrindEQ__2_52_} , \eqref{GrindEQ__2_64_}  into the equation \eqref{GrindEQ__2_75_} we get

\begin{equation} \label{GrindEQ__2_76_}
L_{\overline{C}} =\sinh \theta L_{V_{1} } -\cosh \theta L_{V_{3} } +\theta ^{*} \left(-\cosh \theta \lambda _{V_{1} } +\sinh \theta \lambda _{V_{3} } \right)
\end{equation}
If we use the equations \eqref{GrindEQ__2_57_} and \eqref{GrindEQ__2_69_} into the equation \eqref{GrindEQ__2_76_} and necessary operations have been done, we get

\begin{multline} \label{GrindEQ__2_77_}
L_{\overline{C}} =\sinh \left(\theta +\varphi \right)L_{U_{1} } +\cosh \left(\theta +\varphi \right)L_{U_{3} } -\\
-\left(\varphi ^{*} +\theta ^{*} \right)\left(\cosh \left(\theta +\varphi \right)\lambda _{U_{1} } +\sinh \left(\theta +\varphi \right)\lambda _{U_{3} } \right)
\end{multline}
For the dual angle of the pitch of the closed ruled surface  , we may write

\[\Lambda _{\overline{C}} =-\left\langle \overrightarrow{D},\overrightarrow{\overline{C}}\right\rangle \]
Because of the equations \eqref{GrindEQ__2_21_} and \eqref{GrindEQ__2_71_} we can obtain\textbf{  }

\begin{equation} \label{GrindEQ__2_78_}
\Lambda _{\overline{C}} =-\sinh \Theta \, \oint Qdt +\cosh \Theta \, \oint Pdt
\end{equation}
If we use the equations \eqref{GrindEQ__2_55_} and \eqref{GrindEQ__2_67_}   into the last equations, we get

\begin{equation} \label{GrindEQ__2_79_}
\Lambda _{\overline{C}} =\sinh \Theta \, \Lambda _{V_{1} } -\cosh \Theta \, \Lambda _{V_{3} }
\end{equation}
If we use the equations \eqref{GrindEQ__2_56_} and \eqref{GrindEQ__2_68_}   into the equations \eqref{GrindEQ__2_79_}, we get

\begin{equation} \label{GrindEQ__2_80_}
\Lambda _{\overline{C}} =\sinh \left(\Theta +\Phi \right)\Lambda _{U_{1} } \, +\cosh \left(\Theta +\Phi \right)\, \Lambda _{U_{3} }
\end{equation}
For the drall of the closed surface , we may write

\[P_{\overline{C}} =\frac{\left\langle d\overrightarrow{\overline{c}},d\overrightarrow{\overline{c}}^{*} \right\rangle }{\left\langle d\overrightarrow{\overline{c}},d\overrightarrow{\overline{c}}\right\rangle } \]

\begin{equation} \label{GrindEQ__2_81_}
P_{\overline{C}} =\frac{-\theta '\theta ^{*'} +\left(p\sinh \theta -q\cosh \theta \right)\left[\left(p^{*} -q\theta ^{*} \right)\sinh \theta +\left(p\theta ^{*} -q^{*} \right)\cosh \theta \right]}{\left(p\sinh \theta -q\cosh \theta \right)^{2} -\theta '^{2} }
\end{equation}
\textbf{b)}If the instantaneous dual Pfaffion vector  $\overrightarrow{\overline{\Psi }}$ is timelike ($\left|P\right|<\left|Q\right|$)

\[P=\left\| \overrightarrow{\overline{\Psi }}\right\| \sinh \Theta \, \, \, \, \, \, ,\, \, \, \, \, \, \, Q=\left\| \overrightarrow{\overline{\Psi }}\right\| \cosh \Theta \]
On the way $\overrightarrow{\overline{C}}=\overrightarrow{\overline{c}}+\varepsilon \overrightarrow{\overline{c}}^{*} $, unit vector  about the vector $\overrightarrow{\overline{\Psi }}$ direction is

\begin{equation} \label{GrindEQ__2_82_}
\overrightarrow{\overline{C}}=\cosh \Theta \, \overrightarrow{V_{1} }-\sinh \Theta \overrightarrow{V_{3} }
\end{equation}
Setting by the values of the vectors  $\overrightarrow{V_{1} }{\rm \; and\; }\overrightarrow{V_{3} }$ as the equations \eqref{GrindEQ__2_38_}  into the equation \eqref{GrindEQ__2_82_}, we get

\[\overrightarrow{\overline{C}}=\left(\cosh \Theta \cosh \Phi +\sinh \Theta \sinh \Phi \right)\overrightarrow{U_{1} }+\left(\sinh \Theta \cosh \Phi +\cosh \Theta \sinh \Phi \right)\overrightarrow{U_{3} }\]

\begin{equation} \label{GrindEQ__2_83_}
\overrightarrow{\overline{C}}=\cosh \left(\Theta +\Phi \right)\overrightarrow{U_{1} }+\sinh \left(\Theta +\Phi \right)\overrightarrow{U_{3} }
\end{equation}
If the equation \eqref{GrindEQ__2_82_} is separated into the dual and real part, we can obtain

\begin{equation} \label{GrindEQ__2_84_}
\left\{\begin{array}{l} {\overrightarrow{\overline{c}}=\cosh \theta \overrightarrow{v_{1} }-\sinh \theta \overrightarrow{v_{3} }} \\ {\overrightarrow{\overline{c}}^{*} =\cosh \theta \overrightarrow{v_{1} }^{*} -\sinh \theta \overrightarrow{v_{3} }^{*} +\theta ^{*} \sinh \theta \overrightarrow{v_{1} }-\theta ^{*} \cosh \theta \overrightarrow{v_{3} }} \end{array}\right.
\end{equation}
The pitch of the closed surface  is obtained as

\[L_{\overline{C}} =\left\langle \overrightarrow{d},\overrightarrow{\overline{c}}^{*} \right\rangle +\left\langle \overrightarrow{d}^{*} ,\overrightarrow{\overline{c}}\right\rangle \]
Setting by the values of the statements $\overrightarrow{d}{\rm \; and\; }\overrightarrow{d}^{*} $ as the equations \eqref{GrindEQ__2_51_}  into the last equations and if we do the necessary operations , we get

\begin{equation} \label{GrindEQ__2_85_}
L_{\overline{C}} =-\sinh \theta \oint p^{*} dt +\cosh \theta \oint q^{*} dt +\theta ^{*} \left(\sinh \theta \oint qdt- \cosh \theta \oint pdt \right)
\end{equation}
or

\begin{equation} \label{GrindEQ__2_86_}
L_{\overline{C}} =\cosh \theta L_{V_{1} } -\sinh \theta L_{V_{3} } +\theta ^{*} \left(-\sinh \theta \lambda _{V_{1} } +\cosh \theta \lambda _{V_{3} } \right)
\end{equation}
If we use the equations \eqref{GrindEQ__2_57_} and \eqref{GrindEQ__2_69_} into the equation \eqref{GrindEQ__2_86_} and necessary operations have been done, we get

\begin{multline} \label{GrindEQ__2_87_}
L_{\overline{C}} =\cosh(\theta +\varphi)L_{U_{1}}+\sinh(\theta+ +\varphi )L_{U_{3}}-\\
-({\varphi}'{*}+{\theta}^{*})(\sinh(\theta +\varphi)\lambda_{U_{1}}+\cosh(\theta +\varphi)\lambda_{U_{3}})
\end{multline}
For the dual angle of the pitch of the closed ruled surface  , we may write

\[\Lambda _{\overline{C}} =-\left\langle \overrightarrow{D},\overrightarrow{\overline{C}}\right\rangle \]
Because of the equations \eqref{GrindEQ__2_50_} and \eqref{GrindEQ__2_82_} we can obtain\textbf{  }

\begin{equation} \label{GrindEQ__2_88_}
\Lambda _{\overline{C}} =-\cosh \Theta \, \oint Qdt +\sinh \Theta \, \oint Pdt
\end{equation}
If we use the equations \eqref{GrindEQ__2_55_} and \eqref{GrindEQ__2_67_}   into the last equations, we get

\begin{equation} \label{GrindEQ__2_89_}
\Lambda _{\overline{C}} =\cosh \Theta \, \Lambda _{V_{1} } -\sinh \Theta \, \Lambda _{V_{3} }
\end{equation}
If we use the equations \eqref{GrindEQ__2_57_} and \eqref{GrindEQ__2_68_}   into the equations \eqref{GrindEQ__2_89_}, we get

\begin{equation} \label{GrindEQ__2_90_}
\Lambda _{\overline{C}} =\cosh \left(\Theta +\Phi \right)\Lambda _{U_{1} } \, +\sinh \left(\Theta +\Phi \right)\, \Lambda _{U_{3} }
\end{equation}
For the drall of the closed surface , we may write

\[P_{\overline{C}} =\frac{\left\langle d\overrightarrow{\overline{c}},d\overrightarrow{\overline{c}}^{*} \right\rangle }{\left\langle d\overrightarrow{\overline{c}},d\overrightarrow{\overline{c}}\right\rangle } \]

\begin{equation} \label{GrindEQ__2_91_}
P_{\overline{C}} =\frac{\theta '\theta ^{*'} +\left(p\cosh \theta -q\sinh \theta \right)\left[\left(p\theta ^{*} -q^{*} \right)\sinh \theta +\left(p^{*} -q\theta ^{*} \right)\cosh \theta \right]}{\left(p\cosh \theta -q\sinh \theta \right)^{2} +\theta '^{2} }
\end{equation}

\bibliographystyle{amsalpha}
    
\end{document}